\documentclass[12pt]{amsart}
\usepackage{amsbsy,amsfonts}
\usepackage{latexsym,eucal,amsmath,amsthm}
\usepackage{amssymb}
\newcommand{\db}{\mathbb }

\def\Z{{\db Z}}

\def\R{{\db R}}
\def\C{{\db C}}
\def\torus{\mathbb T}
\def\integers{\mathbb Z}

\def\reals{\mathbb R}

\def\rp{{^{-1}}}

\def\eps{\varepsilon}
\def\be#1{\begin{equation} \label{#1}}
\def\bas{\begin{align*}}

\newtheorem{theorem}{Theorem}
\newtheorem{proposition}{Proposition}
\newtheorem{lemma}{Lemma}

\theoremstyle{definition}

\theoremstyle{remark}
\newtheorem{remark}{Remark}

\numberwithin{equation}{section}
\numberwithin{lemma}{section}
\numberwithin{remark}{section}

\begin{document}

\title[Instability of the nonlinear Schr\"odinger equation]
{Instability of the periodic \\ nonlinear Schr\"odinger equation}

\author{Michael Christ}
\thanks{M.C.\ is supported in part by N.S.F. grant DMS 9970660.}
\address{University of California, Berkeley}
\author{James Colliander}
\thanks{J.C.\ is supported in part by N.S.F. grant DMS 0100595 and N.S.E.R.C.
grant RGPIN 250233-03.}
\address{University of Toronto}
\author{Terence Tao}
\thanks{T.T.\ is a Clay Prize Fellow and is supported in part by a grant
from the Packard Foundation.}
\address{University of California, Los Angeles}


\begin{abstract}  We study the periodic non-linear Schrodinger equation $-iu_t + u_{xx} = \pm |u|^{p-1} u$
for initial data which are assumed to be small in some negative order  Sobolev 
space $H^s(\torus)$ ($s < 0$), but which may have large $L^2$ mass.  In \cite{burqetal}, \cite{cct1} these equations were
shown to be ill-posed in $H^s(\torus)$, in the sense that the solution map was not uniformly continuous from $H^s(\torus)$
to itself even for short times and small norms.  Here we show that these equations are even more unstable, in different ways for different $p$.
For instance,
for the cubic equation ($p=3$) we show that the solution map is not continuous 
as a map from $H^s(\torus)$ to 
even the space of distributions $(C^\infty(\torus))^*$.  
For the quintic equation, the solution map fails to be uniformly continuous
from $C^\infty$ to $C^{-\infty}$ in the sense that
there exist pairs of solutions which are uniformly bounded in $H^s$, are arbitrarily
close in the $C^\infty$ topology at time $0$, and fail to be close
in the distribution topology at an arbitrarily small time $t>0$.
\end{abstract}
%the solution map is not uniformly continuous from
%$C^\infty(\torus)$ to the space of distributions $(C^\infty(\torus))^*$. {\bf Terry will let Mike have the final say on how to word this claim, and the 
%corresponding claim in the body of the text.}

\date{September 25, 2003.}

\maketitle

{\small
\begin{quote}
Like unhappy families,
every unfortunate scientific idea is unfortunate in its own way.
\emph{G.~I.~Barenblatt} 
\end{quote}
}
\medskip

\section{Two unstable Cauchy problems for non-linear Schr\"odinger equations}  \label{section:cubic}

This is the third in a series of papers \cite{cct1},\cite{cct2} investigating the failure of 
nonlinear evolution equations (in particular, defocusing non-linear Schr\"odinger equations)
to be locally well-posed in certain Sobolev spaces $H^s$.  Informally, we say that a Cauchy problem
is locally well-posed in $H^s$ if for any choice of initial data $u_0$ in $H^s$, there exists a time 
$T = T(\|u_0\|_{H^s}) > 0$ depending only on the norm of the initial data, for which a solution exists
on the time interval $[0,T]$,
is unique, and depends continuously on the initial data as a map from $H^s_x$ to $C^0_t H^s_x$ on
the time interval $[0,T]$,
and we say the Cauchy problem is ill-posed if it is not well-posed.  At first glance, this seems like
a simple enough classificiation of Cauchy problems into two distinct classes.  However, it appears that
there are in fact many different types of ill-posedness, caused by different mechanisms of the
underlying equation (blowup, instability of soliton-type solutions, energy transfer from low 
frequencies to high, energy transfer from high frequencies to low, 
and rapid phase decoherence, to name just a few).   
Even within a single class of equations, such
as the semi-linear defocusing Schr\"odinger equations, the nature of the ill-posedness and well-posedness
seems to vary quite substantially depending on the regularity $s$, the dimension $d$, and the power $p$
of the non-linearity, and it seems of interest not only to determine
the exact ranges of these parameters
for which well-posedness or ill-posedness occurs, but also to understand the underlying mechanisms which
generate ill-posedness, or ensure well-posedness.

In this paper we focus attention on one-dimensional periodic non-linear Schr\"odinger equations.
More specifically, we consider the Cauchy problem
\begin{equation} \label{NLS} \tag{NLS}
\left\{
\begin{aligned}
-iu_t + u_{xx} &= \omega |u|^{p-1}u
\\
u(0,x)&=f(x)\in H^s(\torus)
\end{aligned}
\right.
\end{equation}
where $\torus = \reals^1/2\pi \integers$, $u(t,x)$ is a function from $\R \times \torus$ (or $[0,T] \times \torus$)
to $\C$, $s$ is a real number, and $H^s(\torus)$ is the usual Sobolev space with norm
$$ \| f \|_{H^s(\torus)} := (\sum_{n \in \Z} (1 + |n|)^{2s} |\hat f(n)|^2)^{1/2}$$
with
$$ \hat f(n) := \frac{1}{2\pi} \int_{\torus} f(x) e^{-inx}\ dx.$$
We restrict the parameter $\omega$ to either $\omega = +1$ (defocusing case) or
$\omega = -1$ (focusing case), and restrict the exponent $p$ to odd
integers with the $p = 3$ (cubic nonlinearity) and $p = 5$ (quintic
nonlinearity) described in detail.

We begin our discussion with the cubic equation $p=3$, with either choice\footnote{Heuristically, the choice of sign $\omega$
should not be relevant because in the Hamiltonian $\int_\torus \frac{1}{2} |u_x|^2 + \frac{\omega}{4} |u|^4\ dx$,
one can control the potential energy term by the kinetic term (and the mass $\int_\torus |u|^2$) via the
Gagliardo-Nirenberg inequality.  The choice of sign is however somewhat important in the inverse
scattering theory; see \cite{ablowitz}.} of sign $\omega = \pm 1$.  This equation is 
somewhat special, being completely integrable (see e.g. \cite{ablowitz})
and enjoying a large number of conserved quantities.  In particular, it is possible to prove a bound of the form
$\sup_{t \in \R} \| u(t) \|_{H^k(\torus)} \leq C(k, \|f \|_{H^k(\torus)})$
for all $k=0,1,2\ldots$ and all $H^k$ solutions to \eqref{NLS}.  In light of this, it is perhaps not surprising
that global well-posedness in $H^s(\torus)$ has been established\footnote{It should also be mentioned
that there are certain invariant Gibbs-type measures known for the cubic NLS, see for instance \cite{bourgain:invariant}
for a construction of a measure whose support consists almost surely of functions which are roughly of the Sobolev
regularity of $H^{1/2}(\torus)$.  While the notion of an invariant measure is in some sense a very similar notion to
global well-posedness, the connection between the two still seems to be rather poorly understood.}
 \cite{bourgain} for all $s \geq 0$, although
the methods used in \cite{bourgain} are Fourier-analytic and iterative and do not directly exploit the
complete integrability.  In fact, the result in \cite{bourgain} shows that the solution map is not only continuous
from $H^s(\torus)$ to $C^0_t H^s(\torus)$, but is Lipschitz continuous and even real analytic on bounded
subsets of $H^s(\torus)$.  In particular, it is certainly uniformly continuous on bounded subsets of $H^s(\torus)$.

On the other hand, for any exponent $s<0$, this problem is ill-posed in the Sobolev space $H^s(\torus)$ \cite{burqetal}
(see also \cite{cct1})
in the sense that the mapping from  datum $u(0,\cdot)\in C^\infty(\torus)$
to solution $u(t,\cdot)\in C^\infty(\torus)$ 
fails to be uniformly continuous as a mapping from $H^s(\torus)$ to $H^s(\torus)$; this ill-posedness result is
also related to the Gallilean invariance of \eqref{NLS}.
More precisely, for any $\rho>0$ there exists $\eps>0$ such that for any $\delta>0$
there exist initial data $f_1$ and $f_2\in C^\infty$
satisfying $\|f_j\|_{H^s}\le\rho$ and
$\|f_1-f_2\|_{H^s}<\delta$,
a time $T<\delta$, and corresponding solutions $u_j$
such that $\|u_1(T)-u_2(T)\|_{H^s}>\eps$.
It has also been shown (for the Cauchy problem on $\reals\times\reals$
instead of on $\reals\times\torus$; see \cite{cct1} or \cite{cct2}) that for $s\le -\tfrac12$, the
solution map exhibits rapid {\em norm inflation}:
initial data arbitrarily small in $H^s$  norm can give rise to solutions arbitrarily
large in $H^s$ after an arbitrarily short time.  Observe that this is a more dramatic failure
of ill-posedness than failure of uniform continuity, and in particular implies that the solution map
is discontinuous in $H^s$ at the origin $f = 0$.  That result is based upon a transfer of 
energy from high-frequency spatial modes to low-frequency modes.

In this paper we establish an even more extreme form of discontinuity (without growth
of norms) in the regime $s < 0$,
which roughly
speaking asserts that the solution map is not continuous from $H^s(\torus)$ to even the space of distributions
$(C^\infty(\torus))^*$.

\begin{theorem} \label{thm:cubic}
Let $s<0$, $\omega=\pm 1$, and $p \geq 3$ be an odd integer.
Then
for any $\rho>0$ there exists
a solution $u$ of \eqref{NLS}
satisfying $\|u(0)\|_{H^s(\torus)}\le\rho$,
such that for any $\delta>0$
there exists a solution $\tilde u$
satisfying
\[\|u(0)-\tilde u(0)\|_{H^s(\torus)}<\delta,\]
such that\footnote{Here and in the sequel, $c$ and $C$ denote various positive constants depending only on $s$.
Typically $C$ is large, and $c$ is small.}
\[
\sup_{0 \leq t \leq \delta} |\widehat{u(t)}(0)-\widehat{\tilde u(t)}(0)| >c \rho.\]
\end{theorem}

In particular, the solution map fails to be continuous as a map from
$H^s(\torus)$ to $H^\sigma(\torus)$, no matter how close to $-\infty$ the exponent $\sigma$
may be.
The underlying principle here is that a 
high-frequency mode can drive 
oscillation of a low-frequency modes at a rate highly sensitive to
the amplitude of the high-frequency mode.

We should caution that while we establish that the solution map is highly unstable in the $H^\sigma(\torus)$
topology, we do not actually establish any blowup or inflation of the $H^\sigma(\torus)$ norm.
In fact, in a future paper we will show
that for $s$ in some range $-\epsilon_0 < s < 0$,
it is possible to obtain an \emph{a priori} bound of the form
$$ \sup_{0 \leq t \leq T} \| u(t) \|_{H^s(\torus)} \leq C(s, T, \| u(0) \|_{H^s(\torus)})$$
for all smooth solutions $u$ to \eqref{NLS} with $p=3$.  Thus the (densely defined) solution map from 
$H^s(\torus)$ to $H^s(\torus)$ is bounded, even though it is discontinuous.

Theorem~\ref{thm:cubic} implies instability from a fixed Sobolev space $H^s(\torus)$ to a Sobolev 
space $H^\sigma(\torus)$ with large negative exponent $\sigma$.  This should be contrasted
with the results of Lebeau \cite{lebeau1},\cite{lebeau2},
who (in the context of the supercritical wave equation) proved a roughly dual version of this instability,
that the solution map was unstable from high order Sobolev spaces to the energy space.
Specifically, Lebeau showed that for the nonlinear wave equation $\square u +u^7=0$ in $\reals^{1+3}$,
there exist weak solutions $u_j$ whose initial data are 
$O(1)$ in the energy norm
$E^2(f) = \int_{\reals^3} \tfrac12|\partial_t f|^2 + \tfrac12|\nabla_x f|^2 + \tfrac18|f|^8$,
and differ by an arbitrarily small amount in $H^s$ for any large fixed $s$,
for which  $u_1-u_2$ is $\gtrsim 1$ 
in the energy norm after an arbitrarily short time.
We emphasize however that Lebeau's theorem seems to be more subtle, relying on nonlinear geometric optics; our arguments exploit
the dispersive character of the linear Schr\"odinger equation
and fail to apply directly to wave-type equations.

We now turn our attention to the case of the quintic NLS ($p=5$), again with either choice of sign
$\omega = \pm 1$.  This equation is not believed to be completely integrable (for instance, it fails the
Painlev\'e test for integrability, see e.g. \cite{as}).  For this 
equation the $L^2$ norm (which is conserved) is scale invariant, and 
so one would naturally conjecture that local and global well-posedness hold in $H^s(\torus)$ for all $s \geq 0$.  At
the level of $L^2$ this is not known, even for small initial data (mainly because of the ``logarithmic'' failure 
of a key periodic $L^6_{t,x}$ Strichartz estimate, see \cite{bourgain}), however local well-posedness is known in $H^s(\torus)$
for all $s > 0$, see \cite{bourgain}.  Global well-posedness has also
been established \cite{ckstt} in $H^s(\torus)$ for $s>2/3$ and
recently \cite{bourgain:imethod} for all 
$s > s_*$ for some $s_* < 1/2$; for $s \geq 1$ this follows from the local well-posedness
result of \cite{bourgain} and conservation of the Hamiltonian.

When $s < 0$ it is immediate from the arguments given for the cubic Schr\"odinger equation in \cite{burqetal}
or \cite{cct1} that the solution map is not uniformly continuous from $H^s(\torus)$ to $H^s(\torus)$.
Theorem \ref{thm:cubic} also applies to the quintic equation.  
However for this equation we can also prove 
a stronger form of instability, which roughly speaking
asserts that 
solutions which are $O(1)$ in $H^s$ and are arbitrarily close together in
the $C^\infty$ topology at time $t=0$ do not remain close, even in the distribution
topology, for arbitrarily small times $t>0$.
A precise statement is as follows.
 
%if the only bound imposed on the initial data is an $H^s(\torus)$ bound, then the solution map is not 
%uniformly continuous from $C^{\infty}(\torus)$ to the space of distributions $(C^\infty(\torus))^*$, even for extremely short times.

\begin{theorem} \label{thm:quintic}
Let $\omega=\pm 1$, $p \geq 5$ be an odd integer, and $s<0$. For any $\rho >0$ and $\delta > 0$ there exist 
smooth solutions $u_1,u_2$ of \eqref{NLS} such that $u_1(0) - u_2(0)$ is equal to a constant of magnitude
at most $\delta$, and
\begin{align*}
\|u_1(0)\|_{H^s}+\|u_2(0)\|_{H^s}\le\rho,
\\
\sup_{0 \leq t \leq \delta} |\widehat{u_1(t)}(0)-\widehat{u_2(t)}(0)|\geq c \rho.
\end{align*}
\end{theorem}

The conclusion of Theorem~\ref{thm:quintic} is stronger than that of Theorem \ref{thm:cubic} in that
one has instability even when the initial data are close to each other in a smooth norm, but
is weaker in that we compare two nearby initial 
data $u_1,u_2$ (both depending on $\delta$), rather than comparing one reference datum $u$ independent of $\rho$ 
to a sequence of data $\tilde u$ approaching $u$.  In other words the former result disproves uniform stability, 
whereas the latter disproves pointwise stability.   The mechanism
for instability is slightly different in the two cases; for the cubic equation we have a small (in $H^s(\torus)$ norm) 
amplitude
oscillation at high frequency driving a significant fluctuation at low frequency, while for the quintic equation 
we have a small amplitude oscillation at low frequency interacting with a small (in $H^s(\torus)$ norm) 
amplitude oscillation at high 
frequency to create a significant fluctuation at low frequencies.  
In both cases, the high frequency component is small in $H^s$ but has a large Fourier coefficient
in absolute terms, which is possible precisely because $s$ is negative.

The two arguments do share a general scheme:
one first constructs an explicit approximate solution which already exhibits the desired
instability behavior, and then uses a (rescaled) energy method to show that the actual solutions stay close
to the approximate solutions.  See also \cite{cct1}, \cite{cct2}.

\section{Instability of the cubic NLS}

In this section we fix $p=3$, $\omega = \pm 1$, $\rho > 0$ and $s < 0$, and analyze certain special
solutions of \eqref{NLS}.  We then use these solutions to prove Theorem \ref{thm:cubic}.
The remaining cases consisting of 
odd powers $p \geq 5$ are discussed in Remark \ref{powers} at the end of the paper.

As was observed by Burq et.\ al.\ \cite{burqetal},
the equation \eqref{NLS} can be solved explicitly if the initial data $f(x) = u_0(t,x)$
is supported on a single frequency mode, thus $f(x) = \alpha e^{iNx}$ for some integer $N$ and some complex
number $\alpha$.  Then  there is the explicit solution
$$ u(t,x) := \alpha e^{iN^2 t + iNx + i \omega |\alpha|^2 t}.$$
These solutions already suffice to demonstrate the failure of uniform continuity of the solution
map in $H^s(\torus)$ for any $s < 0$ (see \cite{burqetal}, \cite{cct1}).  However, they do not manifest
the stronger instability property claimed in Theorem \ref{thm:cubic}.  Instead, we will
consider the next simplest type of initial data, supported on two separate frequency modes.
To begin with we consider data of the form $f(x) = \alpha + \beta e^{ix}$, although we will later
use scaling\footnote{One could also use Gallilean invariance to allow both frequency modes to be
non-zero, although we will not need to do so in this paper.}
arguments to change this to $f(x) = \alpha N + \beta N e^{iNx}$.  One should think of $\alpha,\beta$
as being rather small parameters (with $\alpha$ somewhat smaller than $\beta$), but we will need to analyze this
solution for long times (which will become short times after the rescaling).

We begin by constructing an explicit approximate solution with the above initial data.

\begin{lemma}\label{alpha-beta}  Let $\alpha, \beta$ be complex numbers with $|\alpha|,|\beta| \leq \sigma$
for some $0 < \sigma \ll 1$, and 
let $u_{\alpha,\beta}: \R \times \torus
\to \C$ be the function
\begin{equation}\label{u-def}
 u_{\alpha,\beta}(t,x) = \alpha e^{i\omega(|\alpha|^2 + 2|\beta|^2)t} + \beta e^{i\omega(2|\alpha|^2 + |\beta|^2)t} 
e^{ix + it}.
\end{equation}
Then the solution $u'_{\alpha,\beta}$ to the inhomogeneous linear Schr\"odinger equation
\begin{equation*}
\left\{
\begin{aligned}
(-i\partial_t + \partial_{xx})u'_{\alpha,\beta} &= \omega |u_{\alpha,\beta}|^2 u_{\alpha,\beta}\\
u'_{\alpha,\beta}(0,x) &= \alpha + \beta e^{ix}
\end{aligned}
\right.
\end{equation*}
satisfies the long-time estimate
\begin{equation}\label{long-time}
 \sup_{t \in \R} \| u'_{\alpha,\beta}(t) - u_{\alpha,\beta}(t) \|_{H^1(\torus)} \leq C \sigma^3.
\end{equation}
\end{lemma}
\noindent $u_{\alpha,\beta}$ is in this sense an approximate solution to \eqref{NLS} with $p=3$ and 
initial data $\alpha + \beta e^{ix}$.  For an explanation of the non-linear phase factors in \eqref{u-def}, see 
Remark \ref{powers} at the end of this paper.

\begin{proof}
A direct computation shows that
$$ (-i\partial_t + \partial_{xx}) u_{\alpha,\beta}
= \omega (|\alpha|^2 + 2 |\beta|^2) \alpha e^{i\omega(|\alpha|^2 + 2|\beta|^2)t}
+ \omega (2|\alpha|^2 + |\beta|^2) \beta e^{i\omega(2|\alpha|^2 + |\beta|^2)t} e^{ix + it},$$
while
$$
\omega |u_{\alpha,\beta}|^2 u_{\alpha,\beta} =
\omega (|\alpha|^2 + 2 |\beta|^2) \alpha e^{i\omega(|\alpha|^2 + 2|\beta|^2)t}
+ \omega (2|\alpha|^2 + |\beta|^2) \beta e^{i\omega(2|\alpha|^2 + |\beta|^2)t} e^{ix + it} + \omega E$$
where the error term $E$ is equal to
$$
E(t,x) := \alpha^2 \overline{\beta} e^{3i\omega |\beta|^2 t} e^{-ix-it} + 
\beta^2 \overline{\alpha} e^{3i\omega |\alpha|^2 t} e^{2ix + 2it}.
$$
The quantity $v := u'_{\alpha,\beta} - u_{\alpha,\beta}$ thus obeys the inhomogeneous equation
$$
(-i\partial_t + \partial_{xx})v = \omega E; \quad v(0,x) = 0$$
and thus has an explicit solution
$$ v(t,x) = \omega \alpha^2 \overline{\beta} \frac{e^{3i\omega |\beta|^2 t} e^{-2it} - 1}{
-2 + 3\omega |\beta|^2} e^{-ix+it}
+ \omega \beta^2 \overline{\alpha} \frac{e^{3i\omega |\alpha|^2 t} e^{-2it} - 1}{-2 + 3\omega|\alpha|^2}
e^{2ix+4it}.$$
The claim follows since $\sigma$ was assumed to be small.  (In the focusing case $\omega = -1$
the smallness assumption is unnecessary).
\end{proof}

We now use the energy method in a standard manner
to convert this approximate solution to an exact solution. 
Recall that \eqref{NLS} with $p=3$ is globally well-posed in $H^1(\torus)$; this guarantees
that the solution $U_{\alpha,\beta}$ in the next proposition exists for all time and
varies continuously in $H^1$.

\begin{proposition}\label{exact-alphabeta} Let $\alpha, \beta$ be complex numbers with $|\alpha|,|\beta| \leq \sigma$
for some $0 < \sigma \ll 1$, and let $u_{\alpha,\beta}$ be the function defined by \eqref{u-def}.
Let $U_{\alpha,\beta}$ be the exact solution to \eqref{NLS} with $p= 3$,  and initial datum 
$U_{\alpha,\beta}(0,x) = \alpha + \beta e^{ix}$. 
Then 
\begin{equation}\label{sigma-2}
\sup_{0 < t \ll \sigma^{-2} \log(1/\sigma)}  \| U_{\alpha,\beta}(t) - u_{\alpha,\beta}(t) \|_{H^1(\torus)}
\leq C \sigma^3.
\end{equation}
\end{proposition}
\noindent
The notation  $0<t\ll\sigma^{-2}\log(1/\sigma)$ means that $0<t\le c\sigma^{-2}\log(1/\sigma)$
for some sufficiently small constant $c$.

\begin{proof}
As observed above, $U_{\alpha,\beta}$ is globally defined and lies in $H^1(\torus)$ for all times $t$.  
Let $w$ denote the function
$$ w := U_{\alpha,\beta} - u'_{\alpha,\beta}$$
where $u'_{\alpha,\beta}$ was defined in the previous lemma; then $w(t)$ also lies in $H^1(\torus)$ for all
$t$, and obeys the difference equation
\begin{align*}
(-i\partial_t + \partial_{xx}) w &= \omega (|U_{\alpha,\beta}|^2 U_{\alpha,\beta} - |u_{\alpha,\beta}|^2 u_{\alpha,\beta})\\
w(0,x) &= 0.
\end{align*}
In particular, we observe that $\|w(0)\|_{H^1(\torus)} = 0$, and we have the energy inequality
$$ \partial_t \| w(t) \|_{H^1(\torus)} \leq \| (|U_{\alpha,\beta}|^2 U_{\alpha,\beta} - |u_{\alpha,\beta}|^2 u_{\alpha,\beta})(t) \|_{H^1(\torus)}.$$
Now observe the general inequality
$$ \| |f+g|^2 (f+g) - |f|^2 f \|_{H^1(\torus)} \leq C \|g\|_{H^1(\torus)} (\| f \|_{H^1(\torus)}^2 + \|g\|_{H^1(\torus)}^2);$$
this is basically because $H^1(\torus)$ is closed under multiplication (as follows from the
Leibniz rule and the Sobolev embedding $H^1(\torus) \subset L^\infty(\torus)$).  From this and \eqref{long-time}
we deduce that
$$ \partial_t \| w(t) \|_{H^1(\torus)} \leq 
C (\|w(t) \|_{H^1(\torus)} + \sigma^3) (\| u_{\alpha,\beta}(t) \|_{H^1(\torus)}^2 + (\|w(t) \|_{H^1(\torus)} + \sigma^3)^2).$$
But trivially $\|u_{\alpha,\beta}(t) \|_{H^1(\torus)} \leq C \sigma$, and hence we have
(after some simplifying) the differential inequality
$$ \partial_t \| w(t) \|_{H^1(\torus)} \leq 
C \sigma^5 + C \sigma^2 \| w(t) \|_{H^1(\torus)} + C \|w(t) \|_{H^1(\torus)}^3.$$
The third term on the right-hand side is negligible if $\|w(t) \|_{H^1(\torus)} \ll \sigma$.  Ignoring this term
for the moment, Gronwall's inequality gives
$$ \| w(t) \|_{H^1(\torus)} \leq C \sigma^5 \exp( C \sigma^2 t ) \leq C \sigma^3 \ll \sigma$$
for $0 < t \ll \sigma^{-2} \log(1/\sigma)$; the hypothesis $\|w(t) \|_{H^1(\torus)} \ll \sigma$ can then
be removed by the usual continuity argument.  The estimate \eqref{sigma-2} then follows from this estimate,
\eqref{long-time}, and the triangle inequality.
\end{proof}

We are now ready to prove Theorem \ref{thm:cubic}.

\begin{proof}[Proof of Theorem~\ref{thm:cubic}]
We restrict our attention to the case $p=3$ and describe the
modifications for the other powers later in Remark \ref{powers}.
Fix $\rho$ and $s$.  We may assume that $s$ is close to 0, say $-1/2 < s < 0$, since the Theorem for larger $s$
will clearly imply the same result for lower $s$.
We set $u(t,x)$ to be the explicit solution
$$ u(t,x) := \rho' e^{i \omega |\rho'|^2 t}$$
to \eqref{NLS}, where $\rho' := \frac{1}{4} \rho$.  Then we clearly 
have $\| u(0) \|_{H^s(\torus)} = \frac{1}{4} \rho$.

Now take $\delta > 0$, which we may assume to be much smaller than $\rho$.  We select $\tilde u(t,x)$ to be the solution
to \eqref{NLS} with initial datum
$$ \tilde u(0,x) = \rho' + \frac{1}{4} \delta N^{-s} e^{iNx},$$
where $N$ is a very large parameter depending on $\delta$ and $\rho$ to be chosen later.
Note that $\tilde u(0,x)$ is well-defined thanks to the global well-posedness of \eqref{NLS} in $H^1(\torus)$ (for
instance).  Also it is clear that
$$ \| u(0) \|_{H^s(\torus)} + \| \tilde u(0) \|_{H^s(\torus)} < \rho$$
and
$$ \| u(0) - \tilde u(0) \|_{H^s(\torus)} < \delta$$
if $N$ is large enough.
Now we investigate $\tilde u(t)$ for times $t > 0$.  Define the parameters
$$ \alpha := \frac{\rho'}{N}; \quad \beta := \frac{1}{4N} \delta N^{-s};$$
note that $|\alpha|, |\beta| \leq \sigma$ for some $\sigma = O(\delta N^{-s-1})$ if $N$ is large enough.  Then from 
Proposition \ref{exact-alphabeta} we have the exact solution $U_{\alpha,\beta}$ to \eqref{NLS} with initial
datum $U_{\alpha,\beta}(0,x) = \alpha + \beta e^{ix}$, and this exact solution obeys the bounds \eqref{sigma-2}.
A rescaling argument and the uniqueness of solutions in $H^1(\torus)$ (known from the well-posedness theory)
then shows that
$$ \tilde u(t,x) = N U_{\alpha,\beta}(N^2 t, Nx);$$
the rescaling by $N$ on $\torus$ is justified since $N$ is an integer.  In particular 
$$ \widehat{\tilde u(t)}(0) = N \widehat{U_{\alpha,\beta}(t)}(0).$$
But from \eqref{u-def} 
$$ \widehat{u_{\alpha,\beta}(t)}(0) = \alpha e^{i\omega(|\alpha|^2 + 2|\beta|^2)t}$$
while from \eqref{sigma-2} 
$$ |\widehat{U_{\alpha,\beta}(t)}(0) - \widehat{u_{\alpha,\beta}(t)}(0)| \leq C \sigma^3$$
for all $0 \leq t \ll \sigma^{-2} \log 1/\sigma$.  Meanwhile, 
$$ \widehat{u(t)}(0) = \rho' = N \alpha.$$
Putting this all together, we obtain
$$ |\widehat{u(t)}(0) - \widehat{\tilde u(t)}(0)|
= \rho' |e^{i\omega(|\alpha|^2 + 2|\beta|^2)N^2 t} - e^{i\omega|\alpha|^2N^2 t}| + O(N \sigma^3)
= \rho' |e^{i\omega 2|\beta|^2N^2 t} -1| + O(\delta^3 N^{-3s-2})
$$
for all $0 \leq t \ll N^{-2} \sigma^{-2} \log 1/\sigma$.  Since $2|\beta|^2 \sim \sigma^2$
and $\sigma \sim \delta N^{-s-1} \ll 1$, the phase $e^{i\omega 2|\beta|^2N^2 t}$
undergoes more than one full rotation in the time interval $0 \leq t \ll N^{-2} \sigma^{-2} \log 1/\sigma$,
and so
$$ \sup_{0 \leq t \ll N^{-2} \sigma^{-2} \log 1/\sigma} \rho' |e^{i\omega 2|\beta|^2N^2 t} - 1|
\geq c \rho.$$
But since $s > -1/2$, $N \sigma^3 \sim  \delta^3 N^{-3s-2} = O(N^{-1/2}\delta^3 )\ll \rho$ if $N$ is
large enough, and thus 
$$ \sup_{0 \leq t \ll N^{-2} \sigma^{-2} \log 1/\sigma} 
|\widehat{u(t)}(0) - \widehat{\tilde u(t)}(0)| \geq c\rho.$$
Since the timespan $N^{-2} \sigma^{-2} \log 1/\sigma \lesssim N^{2s}\delta^{-2}(\log N + \log\delta\rp)$ 
tends to zero as $N\to\infty$ while $\delta$ remains fixed,
Theorem \ref{thm:cubic} follows.
\end{proof}

\section{Greater instability for the quintic equation}

We now modify the above analysis for the quintic NLS.  We begin with the analogue of Lemma \ref{alpha-beta}.

\begin{lemma}\label{alpha-beta-quintic}  Let $\alpha, \beta$ be complex numbers with $|\alpha|,|\beta| \leq \sigma$
for some $0 < \sigma \ll 1$, and 
let $u_{\alpha,\beta}: \R \times \torus
\to \C$ be the function
\begin{equation}\label{u-def-quintic}
 u_{\alpha,\beta}(t,x) = \alpha e^{i\omega(|\alpha|^4 + 6 |\alpha|^2 |\beta|^2 + 3|\beta|^4)t} + 
\beta e^{i\omega(|\beta|^4 + 6 |\beta|^2 |\alpha|^2 + 3|\alpha|^4)t} 
e^{ix + it}.
\end{equation}
Then the solution $u'_{\alpha,\beta}$ to the inhomogeneous linear Schr\"odinger equation
\begin{equation*}
\left\{
\begin{aligned}
(-i\partial_t + \partial_{xx}) u'_{\alpha,\beta} &= \omega |u_{\alpha,\beta}|^4 u_{\alpha,\beta}\\
\tilde u_{\alpha,\beta}(0,x) &= \alpha + \beta e^{ix}
\end{aligned}
\right.
\end{equation*}
satisfies
\begin{equation}\label{long-time-quintic}
 \sup_{t \in \R} \| u'_{\alpha,\beta}(t) - u_{\alpha,\beta}(t) \|_{H^1(\torus)} \leq C \sigma^5.
\end{equation}
\end{lemma}

This is proven in exactly the same manner as Lemma \ref{alpha-beta}.  The key computation is that
\begin{align*}
\omega |u_{\alpha,\beta}|^4 u_{\alpha,\beta} =&
\omega (|\alpha|^4 + 6 |\alpha|^2 |\beta|^2 + 3 |\beta|^4) \alpha e^{i\omega(|\alpha|^2 + 2|\beta|^2)t}\\
&+ \omega (|\beta|^4 + 6 |\beta|^2 |\alpha|^2 + 3|\alpha|^4) \beta e^{i\omega(2|\beta|^2 + |\alpha|^2)t} e^{ix + it} 
\\
&+ \omega E
\end{align*}
where the error term $E$ has the form
$$
E(t,x) = \sum_{k = -2,-1,2,3} c_k e^{i\omega \lambda_k t} e^{ik(x+t)}$$
for various explicit quintic expressions $c_k$ of $\alpha$ and $\beta$ of total magnitude $O(\sigma^5)$, and
various explicit quartic expressions $\lambda_k$ of $\alpha$ and $\beta$ of total magnitude $O(\sigma^4)$.
This is best seen by expanding $|u_{\alpha,\beta}|^4 u_{\alpha,\beta}$ using \eqref{u-def-quintic}
and grouping the terms based
on the power of $e^{i(x+t)}$ which appears.  The rest of the proof of Lemma \ref{alpha-beta-quintic}
continues as in Lemma \ref{alpha-beta}; the main point is that the phases $k(x+t)$ in the forcing terms
are very far from being resonant phases $kx + k^2 t$ for $k = -2,-1,2,3$, so there are no linear growth factors
in \eqref{long-time-quintic}.
Details of the proof of the lemma are left to the reader.

By repeating the proof of Proposition \ref{exact-alphabeta} with only minor changes (altering the powers of
$\sigma$ due to the higher power of the nonlinearity), we obtain

\begin{proposition}\label{exact-alphabeta-quintic} Let $\alpha, \beta$ be complex numbers with $|\alpha|,|\beta| \leq \sigma$
for some $0 < \sigma \ll 1$, and let $u_{\alpha,\beta}$ be the function defined by \eqref{u-def-quintic}.
Let $U_{\alpha,\beta}$ be the exact solution to \eqref{NLS} with $p
= 5$, an odd integer, and initial datum 
$U_{\alpha,\beta}(0,x) = \alpha + \beta e^{ix}$. 
Then 
\begin{equation}\label{sigma-4}
\sup_{0 < t \ll \sigma^{-4} \log(1/\sigma)}  \| U_{\alpha,\beta}(t) - u_{\alpha,\beta}(t) \|_{H^1(\torus)}
\leq C \sigma^5.
\end{equation}
\end{proposition}

Quintic NLS is known to be globally well-posed in $H^1(\torus)$ when the conserved $L^2$ norm is small\footnote{In the 
defocusing case $\omega = +1$ one does not need this smallness hypothesis.  See \cite{bourgain}.}, which guarantees existence
of the solution discussed in the next proposition. 

One could now repeat the argument of Theorem \ref{thm:cubic} to prove the obvious analogue of that result
for the quintic equation; see Remark \ref{powers}.  To prove Theorem \ref{thm:quintic} we proceed in a slightly
different manner.

\begin{proof}[Proof of Theorem~\ref{thm:quintic}]
We carry out the analysis in the case $p=5$ and describe the other
cases below in Remark \ref{powers}. 
Fix $\rho$, $\delta$, $\omega$, $s$; again we may assume that $\delta$ is much smaller than $\rho$.
We shall need a large number $M \gg 1$ depending on $\rho$, $\delta$ to be chosen later, and a very
large integer $N \gg 1$ depending on $\rho, \delta, M$ to be chosen later (we will eventually need
$N$ to be exponentially larger than $M$).  
Let $u_1$, $u_2$ be the exact solutions
to \eqref{NLS} with initial data
$$ u_j(0,x) := \rho' + \delta j + \rho M e^{iNx}$$
for $j=1,2$ where $\rho' = \frac{1}{4} \rho$.  
Clearly $u_1(0)$ and $u_2(0)$ differ by the constant $\delta$,
and we have
$$ \| u_1(0) \|_{H^s(\torus)} + \| u_2(0) \|_{H^s(\torus)} \leq \rho$$
if $\delta$ is small enough and $N$ is large enough depending on $M$.

We now rescale as before, although for the quintic equation the scaling is slightly different.  Performing
the scaling computation, we see that
$$  u_j(t,x) = N^{1/2} U_{\alpha_j,\beta}(N^2 t, Nx),$$
is a solution for $j=1,2$, where $U_{\alpha_j,\beta}$ are as in Proposition \ref{exact-alphabeta-quintic} and
$$ \alpha_j := \frac{\rho' + \delta j}{N^{1/2}}; \quad \beta := \frac{\rho M}{N^{1/2}};$$
observe that $|\alpha_j|, |\beta| \leq \sigma$ for some $\sigma \sim \rho M N^{-1/2}$.
As before, 
$$ \widehat{u_j(t)}(0) = N^{1/2} \widehat{U_{\alpha_j,\beta}(t)}(0)$$
while from \eqref{sigma-4} 
$$ | \widehat{U_{\alpha_j,\beta}(t)}(0) - \widehat{u_{\alpha_j,\beta}(t)}(0)| \leq C \sigma^4$$
for all $0 \leq t \ll \sigma^{-4} \log 1/\sigma$.
Finally, from \eqref{u-def-quintic} 
$$ \widehat{u_{\alpha_j,\beta}(t)}(0) = \alpha_j e^{i\omega(|\alpha_j|^4 + 6 |\alpha_j|^2 |\beta|^2 + 3|\beta|^4)t}.$$
Putting this all together, we see that
\begin{multline}\label{u1-u2} | \widehat{u_1(t)}(0) - \widehat{u_2(t)}(0)| 
\\
=
N^{1/2} |\alpha_1 e^{i\omega(|\alpha_1|^4 + 6 |\alpha_1|^2 |\beta|^2 + 3|\beta|^4)N^2 t} -
\alpha_2 e^{i\omega(|\alpha_2|^4 + 6 |\alpha_2|^2 |\beta|^2 + 3|\beta|^4)N^2 t}|
+ O(N^{1/2} \sigma^4)
\end{multline}
for all 
$$
0 \leq t \ll N^{-2} \sigma^{-4} \log 1/\sigma \sim \rho^{-4} M^{-4} \log N.
$$
In particular, if we choose $N$ exponentially large relative to $\rho$ and $M$ then 
\eqref{u1-u2} holds for all $0 \leq t \leq \delta$.
We can simplify the right-hand side of \eqref{u1-u2} somewhat as
\begin{equation}\label{u1-u2-simplified}
 | \widehat{u_1(t)}(0) - \widehat{u_2(t)}(0)| =
N^{1/2} |\alpha_1 e^{i\omega ((|\alpha_1|^4-|\alpha_2|^4) + 6 (|\alpha_1|^2 - |\alpha_2|^2) |\beta|^2)N^2 t} -
\alpha_2|
+ O(\rho^4 M^4 N^{-3/2}).
\end{equation}
Now 
$$ ||\alpha_1|^2 - |\alpha_2|^2| \sim N^{-1} \rho \delta; \quad 
|\alpha_1|^4-|\alpha_2|^4 =O(N^{-2}\rho^3\delta)$$
and hence (for $M$ large enough)
$$ ((|\alpha_1|^4-|\alpha_2|^4) + 6 (|\alpha_1|^2 - |\alpha_2|^2) |\beta|^2)N^2 \sim \rho^3 \delta M^2.$$
In particular, if $M$ is large enough depending on $\rho$, $\delta$ then this phase will complete more than
one full revolution during the time interval $0 \leq t \leq \delta$, and hence
$$ \sup_{0 \leq t \leq \delta} |\alpha_1 e^{i((|\alpha_1|^4-|\alpha_2|^4) + 6 (|\alpha_1|^2 - |\alpha_2|^2) |\beta|^2)N^2 t} -
\alpha_2| = \alpha_1 + \alpha_2 \geq c \rho N^{-1/2}$$
and hence by \eqref{u1-u2-simplified} 
$$ \sup_{0 \leq t \leq \delta} | \widehat{u_1(t)}(0) - \widehat{u_2(t)}(0)| 
\ge c\rho -O(\rho^4M^4 N^{-3/2}) 
\geq c \rho$$
as desired, if $N$ is large enough relative to $M$.
\end{proof}

\begin{remark}
\label{powers}
We briefly describe the modifications to the proofs
given above to cover the cases $p \geq 5$ in Theorem 1 and $p \geq 7$ in Theorem 2. Let $p = 2m
+1$. We recast the (NLS) dynamics in terms of Fourier coefficients as
\begin{equation}
  \label{FNLS}
  \left\{
\begin{aligned}
\partial_t a_k & = i k^2 a_k + i \omega \sum\limits_{k_1 - k_2 + \dots +
  k_{2m+1} = k} a_{k_1} {\overline{a_{k_2}}} \dots a_{k_{2m+1}} \\
a_k (0) & = \widehat{f} (k)
\end{aligned}
\right.
\end{equation}
where
\begin{equation*}
  u(t,x) = \sum_k a_k (t) e^{ikx}.
\end{equation*}
We take initial data satisfying $a_0 (0) = \alpha, ~ a_1 (0) = \beta$
and $a_k (0) = 0$ for all other $k$. Using the intuition that the $k=0$ and $k=1$
Fourier modes should dominate the evolution for this initial data, 
we approximate the \eqref{FNLS}
dynamics by restricting the sum to $k_1, \dots, k_{2m+1} \in \{ 0 ,
1\}$ so that $k_1 - k_2 + \dots + k_{2m+1} \in \{ -m , \dots,
m+1\}$. Thus, there emerges a first order ODE system in $2m+2$
unknowns:
\begin{equation}
  \label{bigode}
  \left\{
\begin{aligned}
\partial_t a_0 & = i \omega \left[ \sum_{j=0}^m 
\binom{m+1}{j}\binom{m}{j}
%{{m+1} \choose j} {m \choose j} 
|a_0|^{2m - 2j} |a_1|^{2j} \right] a_0 \\
\partial_t a_1 & =i a_1 +  i \omega \left[ \sum_{j=0}^m 
\binom{m+1}{j}\binom{m}{j}
%{{m+1} \choose j} {m \choose j} 
|a_1|^{2m - 2j} |a_0|^{2j} \right] a_1 \\
\partial_t a_k & =i a_k +  G_k (a_{-m} , \dots, a_{m+1} ) ~{\mbox{for all
    other}}~k \in \{ -m , \dots , m+1 \}. 
\end{aligned}
\right.
\end{equation}
The equations for $a_0$ and $a_1$ form a self-contained system
whose solution generalizes (and provides motivation for introducing) 
\eqref{u-def} and \eqref{u-def-quintic}; indeed it is easy to verify
the conservation laws $|a_0(t)| = |\alpha|$, $|a_1(t)| = |\beta|$ for this system,
at which point one can explicitly solve for $a_0$ and $a_1$ as
\begin{align*}
a_0(t) & = \alpha \exp\left( i \omega t \sum_{j=0}^m 
\binom{m+1}{j}\binom{m}{j}
%{{m+1} \choose j} {m \choose j} 
|\alpha|^{2m - 2j} |\beta|^{2j}\right)\\
a_1 (t) & = \beta \exp\left( i \omega t \sum_{j=0}^m 
\binom{m+1}{j}\binom{m}{j}
%{{m+1} \choose j} {m \choose j} 
|\beta|^{2m - 2j} |\alpha|^{2j}\right) e^{it}\  .
\end{align*}
Then one could also evaluate $a_k(t)$ explicitly if desired. An
imitation of the discussion in Section 2, with altered powers of
$\sigma$ and the rescaling associated to $p = 2m+1$, proves the
remaining cases of Theorem 1.

For $m \geq 2$, the quintic and higher power cases, the $a_0$ equation
in \eqref{bigode} contains a term of the form $c_{m,j} |a_1|^{2m -2 }
|a_0|^2 a_0.$ Note that such a term does not occur in the cubic case
and this is why no analog of Theorem 2 has been established here. By
choosing, for $j = 1,2$, 
\begin{equation*}
\alpha_j = \frac{\rho ' + \delta j}{N^{1/m}}, ~ \beta = \frac{\rho M}{N^{1/m}},
\end{equation*}
the discussion in Section 3 adapts to prove the remaining cases of
Theorem 2.
\end{remark}

\end{document}